\documentclass{amsart}
\usepackage{amsfonts}

\newtheorem{theorem}{Theorem}
\theoremstyle{plain}

\newtheorem{corollary}{Corollary}

\newtheorem{proposition}{Proposition}
\newtheorem{remark}{Remark}

\numberwithin{equation}{section}
\input{tcilatex}

\begin{document}
\title[Additive Reverses of the Continuous Triangle Inequality]{Additive
Reverses of the Continuous Triangle Inequality for Bochner Integral of
Vector-Valued Functions in Hilbert Spaces}
\author{Sever S. Dragomir}
\address{School of Computer Science and Mathematics\\
Victoria University of Technology\\
PO Box 14428, MCMC 8001\\
Victoria, Australia.}
\email{sever@csm.vu.edu.au}
\urladdr{http://rgmia.vu.edu.au/SSDragomirWeb.html}
\date{April 16, 2004.}

\begin{abstract}
Some additive reverses of the continuous triangle inequality for Bochner
integral of vector-valued functions in Hilbert spaces are given.
Applications for complex-valued functions are provided as well.
\end{abstract}

\maketitle

\section{Introduction}

Let $f:\left[ a,b\right] \rightarrow \mathbb{K}$, $\mathbb{K}=\mathbb{C}$ or 
$\mathbb{R}$ be a Lebesgue integrable function. The following inequality is
the continuous version of the \textit{triangle inequality}%
\begin{equation}
\left\vert \int_{a}^{b}f\left( x\right) dx\right\vert \leq
\int_{a}^{b}\left\vert f\left( x\right) \right\vert dx,  \label{1.1}
\end{equation}%
and plays a fundamental role in Mathematical Analysis and its applications.

It seems, see \cite[p. 492]{MPF}, that the first reverse inequality for (\ref%
{1.1}) was obtained by J. Karamata in his book from 1949, \cite{K}:%
\begin{equation}
\cos \theta \int_{a}^{b}\left\vert f\left( x\right) \right\vert dx\leq
\left\vert \int_{a}^{b}f\left( x\right) dx\right\vert  \label{1.2}
\end{equation}%
provided%
\begin{equation*}
\left\vert \arg f\left( x\right) \right\vert \leq \theta ,\ \ x\in \left[ a,b%
\right] ,
\end{equation*}%
where $\theta $ is a given angle in $\left( 0,\frac{\pi }{2}\right) .$

This integral inequality is the continuous version of a reverse inequality
for the generalised triangle inequality%
\begin{equation}
\cos \theta \sum_{i=1}^{n}\left\vert z_{i}\right\vert \leq \left\vert
\sum_{i=1}^{n}z_{i}\right\vert ,  \label{1.3}
\end{equation}%
provided%
\begin{equation*}
a-\theta \leq \arg \left( z_{i}\right) \leq a+\theta ,\ \ \text{for \ }i\in
\left\{ 1,\dots ,n\right\} ,
\end{equation*}%
where $a\in \mathbb{R}$ and $\theta \in \left( 0,\frac{\pi }{2}\right) ,$
which, as pointed out in \cite[p. 492]{MPF}, was first discovered by M.
Petrovich in 1917, \cite{P}, and, subsequently rediscovered by other
authors, including J. Karamata \cite[p. 300 -- 301]{K}, H.S. Wilf \cite{W},
and in an equivalent form by M. Marden \cite{M}.

The first to consider the problem for sums in the general case of Hilbert
and Banach spaces, were J.B. Diaz and F.T. Metcalf \cite{DM}.

In our previous work \cite{SSD1}, we pointed out some continuous versions of
Diaz and Metcalf reverses of the generalised triangle inequality.

We mention here some results from \cite{SSD1} which may be compared with the
new results obtained in Sections 2 and 3 below.

\begin{theorem}
\label{ta}If $f\in L\left( \left[ a,b\right] ;H\right) ,$ the space of
Bochner integrable functions defined on $\left[ a,b\right] $ and with values
in the Hilbert space $H,$ and there exists a constant $K\geq 1$ and a vector 
$e\in H,$ $\left\Vert e\right\Vert =1$ such that%
\begin{equation}
\left\Vert f\left( t\right) \right\Vert \leq K\func{Re}\left\langle f\left(
t\right) ,e\right\rangle \ \ \ \text{for a.e. }t\in \left[ a,b\right] ,
\label{1.4}
\end{equation}%
then we have the inequality:%
\begin{equation}
\int_{a}^{b}\left\Vert f\left( t\right) \right\Vert dt\leq K\left\Vert
\int_{a}^{b}f\left( t\right) dt\right\Vert .  \label{1.5}
\end{equation}%
The case of equality holds in (\ref{1.5}) if and only if%
\begin{equation}
\int_{a}^{b}f\left( t\right) dt=\frac{1}{K}\left( \int_{a}^{b}\left\Vert
f\left( t\right) \right\Vert dt\right) \cdot e.  \label{1.6}
\end{equation}
\end{theorem}

As particular cases of interest that may be applied in practice, we note the
following corollaries established in \cite{SSD1}.

\begin{corollary}
\label{cb}Let $e$ be a unit vector in the Hilbert space $\left(
H;\left\langle \cdot ,\cdot \right\rangle \right) ,$ $\rho \in \left(
0,1\right) $ and $f\in L\left( \left[ a,b\right] ;H\right) $ so that%
\begin{equation}
\left\Vert f\left( t\right) -e\right\Vert \leq \rho \ \ \ \text{for a.e. }%
t\in \left[ a,b\right] .  \label{1.7}
\end{equation}%
Then we have the inequality%
\begin{equation}
\sqrt{1-\rho ^{2}}\int_{a}^{b}\left\Vert f\left( t\right) \right\Vert dt\leq
\left\Vert \int_{a}^{b}f\left( t\right) dt\right\Vert ,  \label{1.8}
\end{equation}%
with equality if and only if%
\begin{equation}
\int_{a}^{b}f\left( t\right) dt=\sqrt{1-\rho ^{2}}\left(
\int_{a}^{b}\left\Vert f\left( t\right) \right\Vert dt\right) \cdot e.
\label{1.9}
\end{equation}
\end{corollary}

\begin{corollary}
\label{cc}Let $e$ be a unit vector in $H$ and $M\geq m>0.$ If $f\in L\left( %
\left[ a,b\right] ;H\right) $ is such that%
\begin{equation}
\func{Re}\left\langle Me-f\left( t\right) ,f\left( t\right) -me\right\rangle
\geq 0\ \ \ \text{for a.e. }t\in \left[ a,b\right] ,  \label{1.10}
\end{equation}%
or, equivalently,%
\begin{equation}
\left\Vert f\left( t\right) -\frac{M+m}{2}e\right\Vert \leq \frac{1}{2}%
\left( M-m\right) \ \ \ \text{for a.e. }t\in \left[ a,b\right] ,
\label{1.11}
\end{equation}%
then we have the inequality%
\begin{equation}
\frac{2\sqrt{mM}}{M+m}\int_{a}^{b}\left\Vert f\left( t\right) \right\Vert
dt\leq \left\Vert \int_{a}^{b}f\left( t\right) dt\right\Vert ,  \label{1.12}
\end{equation}%
or, equivalently,%
\begin{equation}
0\leq \int_{a}^{b}\left\Vert f\left( t\right) \right\Vert dt-\left\Vert
\int_{a}^{b}f\left( t\right) dt\right\Vert \leq \frac{\left( \sqrt{M}-\sqrt{m%
}\right) ^{2}}{M+m}\left\Vert \int_{a}^{b}f\left( t\right) dt\right\Vert .
\label{1.13}
\end{equation}%
The equality holds in (\ref{2.11}) (or in the second part of (\ref{2.12}) if
and only if%
\begin{equation*}
\int_{a}^{b}f\left( t\right) dt=\frac{2\sqrt{mM}}{M+m}\left(
\int_{a}^{b}\left\Vert f\left( t\right) \right\Vert dt\right) e.
\end{equation*}
\end{corollary}

The main aim of this paper is to point out additive reverses for the
continuous triangle inequality, namely, upper bounds for the nonnegative
difference%
\begin{equation*}
\int_{a}^{b}\left\Vert f\left( t\right) \right\Vert dt-\left\Vert
\int_{a}^{b}f\left( t\right) dt\right\Vert
\end{equation*}%
under various assumptions on the function $f\in L\left( \left[ a,b\right]
;H\right) .$

Both the case for a unit vector $e\in H$ and a family of orthonormal vectors 
$\left\{ e_{i}\right\} _{i\in \left\{ 1,\dots ,n\right\} }$ are analysed.
Applications for complex-valued Lebesgue integrable functions are given as
well.

\section{Some Additive Reverses for a Unit Vector}

The following result holds.

\begin{theorem}
\label{t2.1}If $f\in L\left( \left[ a,b\right] ;H\right) $ is such that
there exists a vector $e\in H,$ $\left\Vert e\right\Vert =1$ and $k:\left[
a,b\right] \rightarrow \lbrack 0,\infty ),$ a Lebesgue integrable function
with%
\begin{equation}
\left\Vert f\left( t\right) \right\Vert -\func{Re}\left\langle f\left(
t\right) ,e\right\rangle \leq k\left( t\right) \ \ \ \text{for a.e. }t\in %
\left[ a,b\right] ,  \label{2.1}
\end{equation}%
then we have the inequality:%
\begin{equation}
\left( 0\leq \right) \int_{a}^{b}\left\Vert f\left( t\right) \right\Vert
dt-\left\Vert \int_{a}^{b}f\left( t\right) dt\right\Vert \leq
\int_{a}^{b}k\left( t\right) dt.  \label{2.2}
\end{equation}%
The equality holds in (\ref{2.2}) if and only if%
\begin{equation}
\int_{a}^{b}\left\Vert f\left( t\right) \right\Vert dt\geq
\int_{a}^{b}k\left( t\right) dt  \label{2.3}
\end{equation}%
and%
\begin{equation}
\int_{a}^{b}f\left( t\right) dt=\left( \int_{a}^{b}\left\Vert f\left(
t\right) \right\Vert dt-\int_{a}^{b}k\left( t\right) dt\right) e.
\label{2.4}
\end{equation}
\end{theorem}

\begin{proof}
If we integrate the inequality (\ref{2.1}), we get%
\begin{equation}
\int_{a}^{b}\left\Vert f\left( t\right) \right\Vert dt\leq \func{Re}%
\left\langle \int_{a}^{b}f\left( t\right) dt,e\right\rangle
+\int_{a}^{b}k\left( t\right) dt.  \label{2.5}
\end{equation}%
By Schwarz's inequality for $e$ and $\int_{a}^{b}f\left( t\right) dt,$ we
have%
\begin{align}
\func{Re}\left\langle \int_{a}^{b}f\left( t\right) dt,e\right\rangle & \leq
\left\vert \func{Re}\left\langle \int_{a}^{b}f\left( t\right)
dt,e\right\rangle \right\vert \leq \left\vert \left\langle
\int_{a}^{b}f\left( t\right) dt,e\right\rangle \right\vert  \label{2.6} \\
& \leq \left\Vert \int_{a}^{b}f\left( t\right) dt\right\Vert \left\Vert
e\right\Vert =\left\Vert \int_{a}^{b}f\left( t\right) dt\right\Vert .  \notag
\end{align}%
Making use of (\ref{2.5}) and (\ref{2.6}), we deduce the desired inequality (%
\ref{2.2}).

If (\ref{2.3}) and (\ref{2.4}) hold true, then%
\begin{equation*}
\left\Vert \int_{a}^{b}f\left( t\right) dt\right\Vert =\left\vert
\int_{a}^{b}\left\Vert f\left( t\right) \right\Vert dt-\int_{a}^{b}k\left(
t\right) dt\right\vert \left\Vert e\right\Vert =\int_{a}^{b}\left\Vert
f\left( t\right) \right\Vert dt-\int_{a}^{b}k\left( t\right) dt
\end{equation*}%
and the equality holds true in (\ref{2.2}).

Conversely, if the equality holds in (\ref{2.2}), then, obviously (\ref{2.3}%
) is valid and we need only to prove (\ref{2.4}).

If $\left\Vert f\left( t\right) \right\Vert -\func{Re}\left\langle f\left(
t\right) ,e\right\rangle <k\left( t\right) $ for a.e. $t\in \left[ a,b\right]
,$ then (\ref{2.5}) holds as a strict inequality, implying that (\ref{2.2})
also holds as a strict inequality. Therefore, if we assume that equality
holds in (\ref{2.2}), then we must have%
\begin{equation}
\left\Vert f\left( t\right) \right\Vert =\func{Re}\left\langle f\left(
t\right) ,e\right\rangle +k\left( t\right) \ \ \text{for a.e.\ \ }t\in \left[
a,b\right] .  \label{2.7}
\end{equation}

It is well known that in Schwarz's inequality $\left\Vert x\right\Vert
\left\Vert y\right\Vert \geq \func{Re}\left\langle x,y\right\rangle $ the
equality holds iff there exists a $\lambda \geq 0$ such that $x=\lambda y.$
Therefore, if we assume that the equality holds in all of (\ref{2.6}), then
there exists a $\lambda \geq 0$ such that%
\begin{equation}
\int_{a}^{b}f\left( t\right) dt=\lambda e.  \label{2.8}
\end{equation}%
Integrating (\ref{2.7}) on $\left[ a,b\right] ,$ we deduce%
\begin{equation*}
\int_{a}^{b}\left\Vert f\left( t\right) \right\Vert dt=\func{Re}\left\langle
\int_{a}^{b}f\left( t\right) dt,e\right\rangle +\int_{a}^{b}k\left( t\right)
dt,
\end{equation*}%
and thus, by (\ref{2.8}), we get%
\begin{equation*}
\int_{a}^{b}\left\Vert f\left( t\right) \right\Vert dt=\lambda \left\Vert
e\right\Vert ^{2}+\int_{a}^{b}k\left( t\right) dt,
\end{equation*}%
giving $\lambda =\int_{a}^{b}\left\Vert f\left( t\right) \right\Vert
dt-\int_{a}^{b}k\left( t\right) dt.$

Using (\ref{2.8}), we deduce (\ref{2.4}) and the theorem is completely
proved.
\end{proof}

The following corollary may be useful for applications.

\begin{corollary}
\label{c2.2}If $f\in L\left( \left[ a,b\right] ;H\right) $ is such that
there exists a vector $e\in H,$ $\left\Vert e\right\Vert =1$ and $\rho \in
\left( 0,1\right) $ such that 
\begin{equation}
\left\Vert f\left( t\right) -e\right\Vert \leq \rho \ \ \ \text{for a.e. }%
t\in \left[ a,b\right] ,  \label{2.9}
\end{equation}%
then we have the inequality%
\begin{align}
& \left( 0\leq \right) \int_{a}^{b}\left\Vert f\left( t\right) \right\Vert
dt-\left\Vert \int_{a}^{b}f\left( t\right) dt\right\Vert  \label{2.10} \\
& \leq \frac{\rho ^{2}}{\sqrt{1-\rho ^{2}}\left( 1+\sqrt{1-\rho ^{2}}\right) 
}\func{Re}\left\langle \int_{a}^{b}f\left( t\right) dt,e\right\rangle  \notag
\\
& \left( \leq \frac{\rho ^{2}}{\sqrt{1-\rho ^{2}}\left( 1+\sqrt{1-\rho ^{2}}%
\right) }\left\Vert \int_{a}^{b}f\left( t\right) dt\right\Vert \right) . 
\notag
\end{align}%
The equality holds in (\ref{2.10}) if and only if%
\begin{equation}
\int_{a}^{b}\left\Vert f\left( t\right) \right\Vert dt\geq \frac{\rho ^{2}}{%
\sqrt{1-\rho ^{2}}\left( 1+\sqrt{1-\rho ^{2}}\right) }\func{Re}\left\langle
\int_{a}^{b}f\left( t\right) dt,e\right\rangle  \label{2.11}
\end{equation}%
and%
\begin{multline}
\int_{a}^{b}f\left( t\right) dt  \label{2.12} \\
=\left( \int_{a}^{b}\left\Vert f\left( t\right) \right\Vert dt-\frac{\rho
^{2}}{\sqrt{1-\rho ^{2}}\left( 1+\sqrt{1-\rho ^{2}}\right) }\func{Re}%
\left\langle \int_{a}^{b}f\left( t\right) dt,e\right\rangle \right) e.
\end{multline}
\end{corollary}

\begin{proof}
Firstly, note that (\ref{2.3}) is equivalent to%
\begin{equation*}
\left\Vert f\left( t\right) \right\Vert ^{2}+1-\rho ^{2}\leq 2\func{Re}%
\left\langle f\left( t\right) ,e\right\rangle ,
\end{equation*}%
giving%
\begin{equation*}
\frac{\left\Vert f\left( t\right) \right\Vert ^{2}}{\sqrt{1-\rho ^{2}}}+%
\sqrt{1-\rho ^{2}}\leq \frac{2\func{Re}\left\langle f\left( t\right)
,e\right\rangle }{\sqrt{1-\rho ^{2}}}
\end{equation*}%
for a.e. $t\in \left[ a,b\right] .$

Since, obviously%
\begin{equation*}
2\left\Vert f\left( t\right) \right\Vert \leq \frac{\left\Vert f\left(
t\right) \right\Vert ^{2}}{\sqrt{1-\rho ^{2}}}+\sqrt{1-\rho ^{2}}
\end{equation*}%
for any $t\in \left[ a,b\right] ,$ then we deduce the inequality%
\begin{equation*}
\left\Vert f\left( t\right) \right\Vert \leq \frac{\func{Re}\left\langle
f\left( t\right) ,e\right\rangle }{\sqrt{1-\rho ^{2}}}\ \ \ \text{for a.e. }%
t\in \left[ a,b\right] ,
\end{equation*}%
which is clearly equivalent to%
\begin{equation*}
\left\Vert f\left( t\right) \right\Vert -\func{Re}\left\langle f\left(
t\right) ,e\right\rangle \leq \frac{\rho ^{2}}{\sqrt{1-\rho ^{2}}\left( 1+%
\sqrt{1-\rho ^{2}}\right) }\func{Re}\left\langle f\left( t\right)
,e\right\rangle
\end{equation*}%
for a.e. $t\in \left[ a,b\right] .$

Applying Theorem \ref{t2.1} for $k\left( t\right) :=\frac{\rho ^{2}}{\sqrt{%
1-\rho ^{2}}\left( 1+\sqrt{1-\rho ^{2}}\right) }\func{Re}\left\langle
f\left( t\right) ,e\right\rangle ,$ we deduce the desired result.
\end{proof}

In the same spirit, we also have the following corollary.

\begin{corollary}
\label{c2.3}If $f\in L\left( \left[ a,b\right] ;H\right) $ is such that
there exists a vector $e\in H,$ $\left\Vert e\right\Vert =1$ and $M\geq m>0$
such that either%
\begin{equation}
\func{Re}\left\langle Me-f\left( t\right) ,f\left( t\right) -me\right\rangle
\geq 0\ \ \ \text{for a.e. }t\in \left[ a,b\right] ,  \label{2.13}
\end{equation}%
or, equivalently,%
\begin{equation}
\left\Vert f\left( t\right) -\frac{M+m}{2}e\right\Vert \leq \frac{1}{2}%
\left( M-m\right) \ \ \ \text{for a.e. }t\in \left[ a,b\right] ,
\label{2.14}
\end{equation}%
then we have the inequality%
\begin{align}
& \left( 0\leq \right) \int_{a}^{b}\left\Vert f\left( t\right) \right\Vert
dt-\left\Vert \int_{a}^{b}f\left( t\right) dt\right\Vert  \label{2.15} \\
& \leq \frac{\left( \sqrt{M}-\sqrt{m}\right) ^{2}}{2\sqrt{mM}}\func{Re}%
\left\langle \int_{a}^{b}f\left( t\right) dt,e\right\rangle  \notag \\
& \left( \leq \frac{\left( \sqrt{M}-\sqrt{m}\right) ^{2}}{2\sqrt{mM}}%
\left\Vert \int_{a}^{b}f\left( t\right) dt\right\Vert \right) .  \notag
\end{align}%
The equality holds in (\ref{2.15}) if and only if%
\begin{equation*}
\int_{a}^{b}\left\Vert f\left( t\right) \right\Vert dt\geq \frac{\left( 
\sqrt{M}-\sqrt{m}\right) ^{2}}{2\sqrt{mM}}\func{Re}\left\langle
\int_{a}^{b}f\left( t\right) dt,e\right\rangle
\end{equation*}%
and%
\begin{equation*}
\int_{a}^{b}f\left( t\right) dt=\left( \int_{a}^{b}\left\Vert f\left(
t\right) \right\Vert dt-\frac{\left( \sqrt{M}-\sqrt{m}\right) ^{2}}{2\sqrt{mM%
}}\func{Re}\left\langle \int_{a}^{b}f\left( t\right) dt,e\right\rangle
\right) e.
\end{equation*}
\end{corollary}

\begin{proof}
The fact that (\ref{2.13}) and (\ref{2.14}) are equivalent is a simple
exercise and we omit the details.

Observe that (\ref{2.13}) is clearly equivalent to%
\begin{equation*}
\left\Vert f\left( t\right) \right\Vert ^{2}+mM\leq \left( M+m\right) \func{%
Re}\left\langle f\left( t\right) ,e\right\rangle
\end{equation*}%
for a.e. $t\in \left[ a,b\right] ,$ giving the inequality%
\begin{equation*}
\frac{\left\Vert f\left( t\right) \right\Vert ^{2}}{\sqrt{mM}}+\sqrt{mM}\leq 
\frac{M+m}{\sqrt{mM}}\func{Re}\left\langle f\left( t\right) ,e\right\rangle
\end{equation*}%
for a.e. $t\in \left[ a,b\right] .$

Since, obviously,%
\begin{equation*}
2\left\Vert f\left( t\right) \right\Vert \leq \frac{\left\Vert f\left(
t\right) \right\Vert ^{2}}{\sqrt{mM}}+\sqrt{mM}
\end{equation*}%
for any $t\in \left[ a,b\right] ,$ hence we deduce the inequality%
\begin{equation*}
\left\Vert f\left( t\right) \right\Vert \leq \frac{M+m}{\sqrt{mM}}\func{Re}%
\left\langle f\left( t\right) ,e\right\rangle \ \ \text{for a.e. }t\in \left[
a,b\right] ,
\end{equation*}%
which is clearly equivalent to%
\begin{equation*}
\left\Vert f\left( t\right) \right\Vert -\func{Re}\left\langle f\left(
t\right) ,e\right\rangle \leq \frac{\left( \sqrt{M}-\sqrt{m}\right) ^{2}}{2%
\sqrt{mM}}\func{Re}\left\langle f\left( t\right) ,e\right\rangle
\end{equation*}%
for a.e. $t\in \left[ a,b\right] .$

Finally, applying Theorem \ref{t2.1}, we obtain the desired result.
\end{proof}

Another direction providing different upper bounds for the positive
difference%
\begin{equation*}
\int_{a}^{b}\left\Vert f\left( t\right) \right\Vert dt-\left\Vert
\int_{a}^{b}f\left( t\right) dt\right\Vert
\end{equation*}%
is outlined in the following.

\begin{corollary}
\label{c2.4}If $f\in L\left( \left[ a,b\right] ;H\right) $ and $r\in
L_{2}\left( \left[ a,b\right] ;H\right) ,$ $e\in H,$ $\left\Vert
e\right\Vert =1$ are such that 
\begin{equation}
\left\Vert f\left( t\right) -e\right\Vert \leq r\left( t\right) \ \ \ \text{%
for a.e. }t\in \left[ a,b\right] ,  \label{2.16}
\end{equation}%
then we have the inequality%
\begin{equation}
\left( 0\leq \right) \int_{a}^{b}\left\Vert f\left( t\right) \right\Vert
dt-\left\Vert \int_{a}^{b}f\left( t\right) dt\right\Vert \leq \frac{1}{2}%
\int_{a}^{b}r^{2}\left( t\right) dt.  \label{2.17}
\end{equation}%
The equality holds in (\ref{2.17}) if and only if%
\begin{equation*}
\int_{a}^{b}\left\Vert f\left( t\right) \right\Vert dt\geq \frac{1}{2}%
\int_{a}^{b}r^{2}\left( t\right) dt
\end{equation*}%
and%
\begin{equation*}
\int_{a}^{b}f\left( t\right) dt=\left( \int_{a}^{b}\left\Vert f\left(
t\right) \right\Vert dt-\frac{1}{2}\int_{a}^{b}r^{2}\left( t\right)
dt\right) e.
\end{equation*}
\end{corollary}

\begin{proof}
The condition (\ref{2.16}) is obviously equivalent to%
\begin{equation*}
\left\Vert f\left( t\right) \right\Vert ^{2}+1\leq 2\func{Re}\left\langle
f\left( t\right) ,e\right\rangle +r^{2}\left( t\right)
\end{equation*}%
for a.e. $t\in \left[ a,b\right] .$

Using the elementary inequality%
\begin{equation*}
2\left\Vert f\left( t\right) \right\Vert \leq \left\Vert f\left( t\right)
\right\Vert ^{2}+1,\ \ t\in \left[ a,b\right] ,
\end{equation*}%
we deduce%
\begin{equation*}
\left\Vert f\left( t\right) \right\Vert -\func{Re}\left\langle f\left(
t\right) ,e\right\rangle \leq \frac{1}{2}r^{2}\left( t\right)
\end{equation*}%
for a.e. $t\in \left[ a,b\right] .$

Applying Theorem \ref{t2.1} for $k\left( t\right) :=\frac{1}{2}r^{2}\left(
t\right) ,$ $t\in \left[ a,b\right] $, we deduce the desired result.
\end{proof}

Finally, we may state and prove the following result as well.

\begin{corollary}
\label{c2.5}If $f\in L\left( \left[ a,b\right] ;H\right) $, $e\in H,$ $%
\left\Vert e\right\Vert =1$ and $M,m:\left[ a,b\right] \rightarrow \lbrack
0,\infty )$ with $M\geq m$ a.e. on $\left[ a,b\right] ,$ are such that $%
\frac{\left( M-m\right) ^{2}}{M+m}\in L\left[ a,b\right] $ and either%
\begin{equation}
\left\Vert f\left( t\right) -\frac{M\left( t\right) +m\left( t\right) }{2}%
e\right\Vert \leq \frac{1}{2}\left[ M\left( t\right) -m\left( t\right) %
\right] \ \ \ \text{for a.e. }t\in \left[ a,b\right] ,  \label{2.18}
\end{equation}%
or, equivalently,%
\begin{equation}
\func{Re}\left\langle M\left( t\right) e-f\left( t\right) ,f\left( t\right)
-m\left( t\right) e\right\rangle \geq 0\ \ \ \text{for a.e. }t\in \left[ a,b%
\right] ,  \label{2.19}
\end{equation}%
hold, then we have the inequality%
\begin{equation}
\left( 0\leq \right) \int_{a}^{b}\left\Vert f\left( t\right) \right\Vert
dt-\left\Vert \int_{a}^{b}f\left( t\right) dt\right\Vert \leq \frac{1}{4}%
\int_{a}^{b}\frac{\left[ M\left( t\right) -m\left( t\right) \right] ^{2}}{%
M\left( t\right) +m\left( t\right) }dt.  \label{2.20}
\end{equation}%
The equality holds in (\ref{2.20}) if and only if%
\begin{equation*}
\int_{a}^{b}\left\Vert f\left( t\right) \right\Vert dt\geq \frac{1}{4}%
\int_{a}^{b}\frac{\left[ M\left( t\right) -m\left( t\right) \right] ^{2}}{%
M\left( t\right) +m\left( t\right) }dt
\end{equation*}%
and%
\begin{equation*}
\int_{a}^{b}f\left( t\right) dt=\left( \int_{a}^{b}\left\Vert f\left(
t\right) \right\Vert dt-\frac{1}{4}\int_{a}^{b}\frac{\left[ M\left( t\right)
-m\left( t\right) \right] ^{2}}{M\left( t\right) +m\left( t\right) }%
dt\right) e.
\end{equation*}
\end{corollary}

\begin{proof}
The condition (\ref{2.18}) is equivalent to%
\begin{multline*}
\left\Vert f\left( t\right) \right\Vert ^{2}+\left( \frac{M\left( t\right)
+m\left( t\right) }{2}\right) ^{2} \\
\leq 2\left( \frac{M\left( t\right) +m\left( t\right) }{2}\right) \func{Re}%
\left\langle f\left( t\right) ,e\right\rangle +\frac{1}{4}\left[ M\left(
t\right) -m\left( t\right) \right] ^{2}
\end{multline*}%
for a.e. $t\in \left[ a,b\right] ,$ and since%
\begin{equation*}
2\left( \frac{M\left( t\right) +m\left( t\right) }{2}\right) \left\Vert
f\left( t\right) \right\Vert \leq \left\Vert f\left( t\right) \right\Vert
^{2}+\left( \frac{M\left( t\right) +m\left( t\right) }{2}\right) ^{2},\ \ \
t\in \left[ a,b\right]
\end{equation*}%
hence%
\begin{equation*}
\left\Vert f\left( t\right) \right\Vert -\func{Re}\left\langle f\left(
t\right) ,e\right\rangle \leq \frac{1}{4}\frac{\left[ M\left( t\right)
-m\left( t\right) \right] ^{2}}{M\left( t\right) +m\left( t\right) }
\end{equation*}%
for a.e. $t\in \left[ a,b\right] .$

Now, applying Theorem \ref{t2.1} for $k\left( t\right) :=\frac{1}{4}\frac{%
\left[ M\left( t\right) -m\left( t\right) \right] ^{2}}{M\left( t\right)
+m\left( t\right) },$ $t\in \left[ a,b\right] $, we deduce the desired
inequality.
\end{proof}

\section{Additive Reverses for orthonormal Families}

We recall that the family of vectors $\left\{ e_{i}\right\} _{i\in \left\{
1,\dots ,n\right\} }$ in the inner product space $\left( H;\left\langle
\cdot ,\cdot \right\rangle \right) $ is orthonormal if%
\begin{equation*}
\left\langle e_{i},e_{j}\right\rangle =0\text{ \ if \ }i\neq j,\ \ i,j\in
\left\{ 1,\dots ,n\right\} \text{ \ and \ }\left\Vert e_{i}\right\Vert =1%
\text{ \ for \ }i\in \left\{ 1,\dots ,n\right\} .
\end{equation*}

The following reverse of the continuous triangle inequality for vector
valued integrals holds.

\begin{theorem}
\label{t3.1}Let $f\in L\left( \left[ a,b\right] ;H\right) ,$ where $H$ is a
Hilbert space over the real or complex number field $\mathbb{K}$, $\left\{
e_{i}\right\} _{i\in \left\{ 1,\dots ,n\right\} }$ an orthonormal family in $%
H$ and $M_{i}\in L\left[ a,b\right] ,$ $i\in \left\{ 1,\dots ,n\right\} .$
If we assume that%
\begin{equation}
\left\Vert f\left( t\right) \right\Vert -\func{Re}\left\langle f\left(
t\right) ,e_{i}\right\rangle \leq M_{i}\left( t\right) \ \ \text{for a.e. }%
t\in \left[ a,b\right] ,  \label{3.1}
\end{equation}%
then we have the inequality%
\begin{equation}
\int_{a}^{b}\left\Vert f\left( t\right) \right\Vert dt\leq \frac{1}{\sqrt{n}}%
\left\Vert \int_{a}^{b}f\left( t\right) dt\right\Vert +\frac{1}{n}%
\sum_{i=1}^{n}\int_{a}^{b}M_{i}\left( t\right) dt.  \label{3.2}
\end{equation}%
The equality holds in (\ref{3.2}) if and only if%
\begin{equation}
\int_{a}^{b}\left\Vert f\left( t\right) \right\Vert dt\geq \frac{1}{n}%
\sum_{i=1}^{n}\int_{a}^{b}M_{i}\left( t\right) dt  \label{3.3}
\end{equation}%
and 
\begin{equation}
\int_{a}^{b}f\left( t\right) dt=\left( \int_{a}^{b}\left\Vert f\left(
t\right) \right\Vert dt-\frac{1}{n}\sum_{i=1}^{n}\int_{a}^{b}M_{i}\left(
t\right) dt\right) \sum_{i=1}^{n}e_{i}.  \label{3.4}
\end{equation}
\end{theorem}

\begin{proof}
If we integrate the inequality (\ref{3.1}) on $\left[ a,b\right] ,$ we get%
\begin{equation*}
\int_{a}^{b}\left\Vert f\left( t\right) \right\Vert dt\leq \func{Re}%
\left\langle \int_{a}^{b}f\left( t\right) dt,e_{i}\right\rangle
+\int_{a}^{b}M_{i}\left( t\right) dt
\end{equation*}%
for each $i\in \left\{ 1,\dots ,n\right\} .$ Summing these inequalities over 
$i$ from $1$ to $n,$ we deduce%
\begin{equation}
\int_{a}^{b}\left\Vert f\left( t\right) \right\Vert dt\leq \frac{1}{n}\func{%
Re}\left\langle \int_{a}^{b}f\left( t\right)
dt,\sum_{i=1}^{n}e_{i}\right\rangle +\frac{1}{n}\sum_{i=1}^{n}%
\int_{a}^{b}M_{i}\left( t\right) dt.  \label{3.5}
\end{equation}%
By Schwarz's inequality for $\int_{a}^{b}f\left( t\right) dt$ and $%
\sum_{i=1}^{n}e_{i},$ we have%
\begin{align}
& \func{Re}\left\langle \int_{a}^{b}f\left( t\right)
dt,\sum_{i=1}^{n}e_{i}\right\rangle  \label{3.6} \\
& \leq \left\vert \func{Re}\left\langle \int_{a}^{b}f\left( t\right)
dt,\sum_{i=1}^{n}e_{i}\right\rangle \right\vert \leq \left\vert \left\langle
\int_{a}^{b}f\left( t\right) dt,\sum_{i=1}^{n}e_{i}\right\rangle \right\vert
\notag \\
& \leq \left\Vert \int_{a}^{b}f\left( t\right) dt\right\Vert \left\Vert
\sum_{i=1}^{n}e_{i}\right\Vert =\sqrt{n}\left\Vert \int_{a}^{b}f\left(
t\right) dt\right\Vert ,  \notag
\end{align}%
since%
\begin{equation*}
\left\Vert \sum_{i=1}^{n}e_{i}\right\Vert =\sqrt{\left\Vert
\sum_{i=1}^{n}e_{i}\right\Vert ^{2}}=\sqrt{\sum_{i=1}^{n}\left\Vert
e_{i}\right\Vert ^{2}}=\sqrt{n}.
\end{equation*}%
Making use of (\ref{3.5}) and (\ref{3.6}), we deduce the desired inequality (%
\ref{3.2}).

If (\ref{3.3}) and (\ref{3.4}) hold, then 
\begin{align*}
\frac{1}{\sqrt{n}}\left\Vert \int_{a}^{b}f\left( t\right) dt\right\Vert & =%
\frac{1}{\sqrt{n}}\left\vert \int_{a}^{b}\left\Vert f\left( t\right)
\right\Vert dt-\frac{1}{n}\sum_{i=1}^{n}\int_{a}^{b}M_{i}\left( t\right)
dt\right\vert \left\Vert \sum_{i=1}^{n}e_{i}\right\Vert \\
& =\left( \int_{a}^{b}\left\Vert f\left( t\right) \right\Vert dt-\frac{1}{n}%
\sum_{i=1}^{n}\int_{a}^{b}M_{i}\left( t\right) dt\right)
\end{align*}%
and the equality in (\ref{3.2}) holds true.

Conversely, if the equality holds in (\ref{3.2}), then, obviously, (\ref{3.3}%
) is valid.

Taking into account the argument presented above for the previous result (%
\ref{3.2}), it is obvious that, if the equality holds in (\ref{3.2}), then
it must hold in (\ref{3.1}) for a.e. $t\in \left[ a,b\right] $ and for each $%
i\in \left\{ 1,\dots ,n\right\} $ and also the equality must hold in any of
the inequalities in (\ref{3.6}).

It is well known that in Schwarz's inequality $\func{Re}\left\langle
u,v\right\rangle \leq \left\Vert u\right\Vert \left\Vert v\right\Vert ,$ the
equality occurs if and only if $u=\lambda v$ with $\lambda \geq 0,$
consequently, the equality holds in all inequalities from (\ref{3.6})
simultaneously iff there exists a $\mu \geq 0$ with%
\begin{equation}
\mu \sum_{i=1}^{n}e_{i}=\int_{a}^{b}f\left( t\right) dt.  \label{3.7}
\end{equation}%
If we integrate the equality in (\ref{3.1}) and sum over $i,$ we deduce%
\begin{equation}
n\int_{a}^{b}f\left( t\right) dt=\func{Re}\left\langle \int_{a}^{b}f\left(
t\right) dt,\sum_{i=1}^{n}e_{i}\right\rangle
+\sum_{i=1}^{n}\int_{a}^{b}M_{i}\left( t\right) dt.  \label{3.8}
\end{equation}%
Replacing $\int_{a}^{b}f\left( t\right) dt$ from (\ref{3.7}) into (\ref{3.8}%
), we deduce%
\begin{align}
n\int_{a}^{b}f\left( t\right) dt& =\mu \left\Vert
\sum_{i=1}^{n}e_{i}\right\Vert ^{2}+\sum_{i=1}^{n}\int_{a}^{b}M_{i}\left(
t\right) dt  \label{3.9} \\
& =\mu n+\sum_{i=1}^{n}\int_{a}^{b}M_{i}\left( t\right) dt.  \notag
\end{align}

Finally, we note that (\ref{3.7}) and (\ref{3.9}) will produce the required
identity (\ref{3.4}), and the proof is complete.
\end{proof}

The following corollaries may be of interest for applications.

\begin{corollary}
\label{c3.2}Let $f\in L\left( \left[ a,b\right] ;H\right) ,$ $\left\{
e_{i}\right\} _{i\in \left\{ 1,\dots ,n\right\} }$ an orthonormal family in $%
H$ and $\rho _{i}\in \left( 0,1\right) ,$ $i\in \left\{ 1,\dots ,n\right\} $
such that%
\begin{equation}
\left\Vert f\left( t\right) -e_{i}\right\Vert \leq \rho _{i}\ \ \text{for
a.e. }t\in \left[ a,b\right] .  \label{3.10}
\end{equation}%
Then we have the inequalities:%
\begin{align}
& \int_{a}^{b}\left\Vert f\left( t\right) \right\Vert dt  \label{3.11} \\
& \leq \frac{1}{\sqrt{n}}\left\Vert \int_{a}^{b}f\left( t\right)
dt\right\Vert +\func{Re}\left\langle \int_{a}^{b}f\left( t\right) dt,\frac{1%
}{n}\sum_{i=1}^{n}\frac{\rho _{i}^{2}}{\sqrt{1-\rho _{i}^{2}}\left( 1+\sqrt{%
1-\rho _{i}^{2}}\right) }e_{i}\right\rangle   \notag \\
& \leq \frac{1}{\sqrt{n}}\left\Vert \int_{a}^{b}f\left( t\right)
dt\right\Vert \left[ 1+\left( \frac{1}{n}\sum_{i=1}^{n}\frac{\rho _{i}^{2}}{%
\sqrt{1-\rho _{i}^{2}}\left( 1+\sqrt{1-\rho _{i}^{2}}\right) }\right) ^{%
\frac{1}{2}}\right] .  \notag
\end{align}%
The equality holds in the first inequality in (\ref{3.11}) if and only if%
\begin{equation*}
\int_{a}^{b}\left\Vert f\left( t\right) \right\Vert dt\geq \func{Re}%
\left\langle \int_{a}^{b}f\left( t\right) dt,\frac{1}{n}\sum_{i=1}^{n}\frac{%
\rho _{i}^{2}}{\sqrt{1-\rho _{i}^{2}}\left( 1+\sqrt{1-\rho _{i}^{2}}\right) }%
e_{i}\right\rangle 
\end{equation*}%
and%
\begin{multline*}
\int_{a}^{b}f\left( t\right) dt \\
=\left( \int_{a}^{b}\left\Vert f\left( t\right) \right\Vert dt-\func{Re}%
\left\langle \int_{a}^{b}f\left( t\right) dt,\frac{1}{n}\sum_{i=1}^{n}\frac{%
\rho _{i}^{2}}{\sqrt{1-\rho _{i}^{2}}\left( 1+\sqrt{1-\rho _{i}^{2}}\right) }%
e_{i}\right\rangle \right) \sum_{i=1}^{n}e_{i}.
\end{multline*}
\end{corollary}

\begin{proof}
As in the proof of Corollary \ref{c2.2}, the assumption (\ref{3.10}) implies%
\begin{equation*}
\left\Vert f\left( t\right) \right\Vert -\func{Re}\left\langle f\left(
t\right) ,e_{i}\right\rangle \leq \frac{\rho _{i}^{2}}{\sqrt{1-\rho _{i}^{2}}%
\left( \sqrt{1-\rho _{i}^{2}}+1\right) }\func{Re}\left\langle f\left(
t\right) ,e_{i}\right\rangle
\end{equation*}%
for a.e. $t\in \left[ a,b\right] $ and for each $i\in \left\{ 1,\dots
,n\right\} .$

Now, if we apply Theorem \ref{t3.1} for 
\begin{equation*}
M_{i}\left( t\right) :=\frac{\rho _{i}^{2}\func{Re}\left\langle f\left(
t\right) ,e_{i}\right\rangle }{\sqrt{1-\rho _{i}^{2}}\left( \sqrt{1-\rho
_{i}^{2}}+1\right) },\ \ i\in \left\{ 1,\dots ,n\right\} ,\ \ t\in \left[ a,b%
\right] ,
\end{equation*}%
we deduce the first inequality in (\ref{3.11}).

By Schwarz's inequality in $H,$ we have%
\begin{align*}
& \func{Re}\left\langle \int_{a}^{b}f\left( t\right) dt,\frac{1}{n}%
\sum_{i=1}^{n}\frac{\rho _{i}^{2}}{\sqrt{1-\rho _{i}^{2}}\left( 1+\sqrt{%
1-\rho _{i}^{2}}\right) }e_{i}\right\rangle \\
& \leq \left\Vert \int_{a}^{b}f\left( t\right) dt\right\Vert \left\Vert 
\frac{1}{n}\sum_{i=1}^{n}\frac{\rho _{i}^{2}}{\sqrt{1-\rho _{i}^{2}}\left( 1+%
\sqrt{1-\rho _{i}^{2}}\right) }e_{i}\right\Vert \\
& =\frac{1}{n}\left\Vert \int_{a}^{b}f\left( t\right) dt\right\Vert \left(
\sum_{i=1}^{n}\left[ \frac{\rho _{i}^{2}}{\sqrt{1-\rho _{i}^{2}}\left( 1+%
\sqrt{1-\rho _{i}^{2}}\right) }\right] ^{2}\right) ^{\frac{1}{2}},
\end{align*}%
which implies the second inequality in (\ref{3.11}).
\end{proof}

The second result is incorporated in:

\begin{corollary}
\label{c3.3}Let $f\in L\left( \left[ a,b\right] ;H\right) ,$ $\left\{
e_{i}\right\} _{i\in \left\{ 1,\dots ,n\right\} }$ an orthonormal family in $%
H$ and $M_{i}\geq m_{i}>0$ such that either%
\begin{equation}
\func{Re}\left\langle M_{i}e_{i}-f\left( t\right) ,f\left( t\right)
-m_{i}e_{i}\right\rangle \geq 0\   \label{3.12a}
\end{equation}%
or, equivalently,%
\begin{equation*}
\left\Vert f\left( t\right) -\frac{M_{i}+m_{i}}{2}\cdot e_{i}\right\Vert
\leq \frac{1}{2}\left( M_{i}-m_{i}\right) \ 
\end{equation*}%
for a.e. $t\in \left[ a,b\right] $ \ and each \ $i\in \left\{ 1,\dots
,n\right\} .$Then we have 
\begin{align}
& \int_{a}^{b}\left\Vert f\left( t\right) \right\Vert dt  \label{3.13} \\
& \leq \frac{1}{\sqrt{n}}\left\Vert \int_{a}^{b}f\left( t\right)
dt\right\Vert +\func{Re}\left\langle \int_{a}^{b}f\left( t\right) dt,\frac{1%
}{n}\sum_{i=1}^{n}\frac{\left( \sqrt{M_{i}}-\sqrt{m_{i}}\right) ^{2}}{2\sqrt{%
m_{i}M_{i}}}e_{i}\right\rangle   \notag \\
& \leq \frac{1}{\sqrt{n}}\left\Vert \int_{a}^{b}f\left( t\right)
dt\right\Vert \left[ 1+\left( \frac{1}{n}\sum_{i=1}^{n}\frac{\left( \sqrt{%
M_{i}}-\sqrt{m_{i}}\right) ^{4}}{4m_{i}M_{i}}\right) ^{\frac{1}{2}}\right] .
\notag
\end{align}%
The equality holds in the first inequality in (\ref{3.13}) if and only if%
\begin{equation*}
\int_{a}^{b}\left\Vert f\left( t\right) \right\Vert dt\geq \func{Re}%
\left\langle \int_{a}^{b}f\left( t\right) dt,\frac{1}{n}\sum_{i=1}^{n}\frac{%
\left( \sqrt{M_{i}}-\sqrt{m_{i}}\right) ^{2}}{2\sqrt{m_{i}M_{i}}}%
e_{i}\right\rangle 
\end{equation*}%
and%
\begin{multline*}
\int_{a}^{b}f\left( t\right) dt \\
=\left( \int_{a}^{b}\left\Vert f\left( t\right) \right\Vert dt-\func{Re}%
\left\langle \int_{a}^{b}f\left( t\right) dt,\frac{1}{n}\sum_{i=1}^{n}\frac{%
\left( \sqrt{M_{i}}-\sqrt{m_{i}}\right) ^{2}}{2\sqrt{m_{i}M_{i}}}%
e_{i}\right\rangle \right) \sum_{i=1}^{n}e_{i}.
\end{multline*}
\end{corollary}

\begin{proof}
As in the proof of Corollary \ref{c2.3}, from (\ref{3.12a}), we have%
\begin{equation*}
\left\Vert f\left( t\right) \right\Vert -\func{Re}\left\langle f\left(
t\right) ,e_{i}\right\rangle \leq \frac{\left( \sqrt{M_{i}}-\sqrt{m_{i}}%
\right) ^{2}}{2\sqrt{m_{i}M_{i}}}\func{Re}\left\langle f\left( t\right)
,e_{i}\right\rangle
\end{equation*}%
for a.e. $t\in \left[ a,b\right] $ and $i\in \left\{ 1,\dots ,n\right\} .$

Applying Theorem \ref{t3.1} for%
\begin{equation*}
M_{i}\left( t\right) :=\frac{\left( \sqrt{M_{i}}-\sqrt{m_{i}}\right) ^{2}}{2%
\sqrt{m_{i}M_{i}}}\func{Re}\left\langle f\left( t\right) ,e_{i}\right\rangle
,\ \ \text{ }t\in \left[ a,b\right] ,\ i\in \left\{ 1,\dots ,n\right\} ,
\end{equation*}%
we deduce the desired result.
\end{proof}

In a different direction, we may state the following result as well.

\begin{corollary}
\label{c3.4}Let $f\in L\left( \left[ a,b\right] ;H\right) ,$ $\left\{
e_{i}\right\} _{i\in \left\{ 1,\dots ,n\right\} }$ an orthonormal family in $%
H$ and $r_{i}\in L^{2}\left( \left[ a,b\right] \right) ,$ $i\in \left\{
1,\dots ,n\right\} $ such that%
\begin{equation*}
\left\Vert f\left( t\right) -e_{i}\right\Vert \leq r_{i}\left( t\right) \ \
\ \text{for a.e. }t\in \left[ a,b\right] \text{ \ and \ }i\in \left\{
1,\dots ,n\right\} .
\end{equation*}%
Then we have the inequality%
\begin{equation}
\int_{a}^{b}\left\Vert f\left( t\right) \right\Vert dt\leq \frac{1}{\sqrt{n}}%
\left\Vert \int_{a}^{b}f\left( t\right) dt\right\Vert +\frac{1}{2n}%
\sum_{i=1}^{n}\left( \int_{a}^{b}r_{i}^{2}\left( t\right) dt\right) .
\label{3.15}
\end{equation}%
The equality holds in (\ref{3.15}) if and only if%
\begin{equation*}
\int_{a}^{b}\left\Vert f\left( t\right) \right\Vert dt\geq \frac{1}{2n}%
\sum_{i=1}^{n}\left( \int_{a}^{b}r_{i}^{2}\left( t\right) dt\right) 
\end{equation*}%
and%
\begin{equation*}
\int_{a}^{b}f\left( t\right) dt=\left[ \int_{a}^{b}\left\Vert f\left(
t\right) \right\Vert dt-\frac{1}{n}\sum_{i=1}^{n}\left(
\int_{a}^{b}r_{i}^{2}\left( t\right) dt\right) \right] \sum_{i=1}^{n}e_{i}.
\end{equation*}
\end{corollary}

\begin{proof}
As in the proof of Corollary \ref{c2.4}, from (\ref{2.16}), we deduce that%
\begin{equation}
\left\Vert f\left( t\right) \right\Vert -\func{Re}\left\langle f\left(
t\right) ,e_{i}\right\rangle \leq \frac{1}{2}r_{i}^{2}\left( t\right)
\label{3.16}
\end{equation}%
for a.e. $t\in \left[ a,b\right] $ and $i\in \left\{ 1,\dots ,n\right\} .$

Applying Theorem \ref{t3.1} for%
\begin{equation*}
M_{i}\left( t\right) :=\frac{1}{2}r_{i}^{2}\left( t\right) ,\ \ t\in \left[
a,b\right] ,\ i\in \left\{ 1,\dots ,n\right\} ,
\end{equation*}%
we get the desired result.
\end{proof}

Finally, the following result holds.

\begin{corollary}
\label{c3.5}Let $f\in L\left( \left[ a,b\right] ;H\right) ,$ $\left\{
e_{i}\right\} _{i\in \left\{ 1,\dots ,n\right\} }$ an orthonormal family in $%
H$, $M_{i},m_{i}:\left[ a,b\right] \rightarrow \lbrack 0,\infty )$ with $%
M_{i}\geq m_{i}$ a.e. on $\left[ a,b\right] $ and $\frac{\left(
M_{i}-m_{i}\right) ^{2}}{M_{i}+m_{i}}\in L\left[ a,b\right] ,$ and either%
\begin{multline}
\left\Vert f\left( t\right) -\frac{M_{i}\left( t\right) +m_{i}\left(
t\right) }{2}e_{i}\right\Vert \leq \frac{1}{2}\left[ M_{i}\left( t\right)
-m_{i}\left( t\right) \right] ^{2}\ \   \label{3.17} \\
\end{multline}%
or, equivalently,%
\begin{equation*}
\func{Re}\left\langle M_{i}\left( t\right) e_{i}-f\left( t\right) ,f\left(
t\right) -m_{i}\left( t\right) e_{i}\right\rangle \geq 0\ \ \ 
\end{equation*}%
for a.e. $t\in \left[ a,b\right] $, \ $i\in \left\{ 1,\dots ,n\right\} $
hold, then we have the inequality%
\begin{equation}
\int_{a}^{b}\left\Vert f\left( t\right) \right\Vert dt\leq \frac{1}{\sqrt{n}}%
\left\Vert \int_{a}^{b}f\left( t\right) dt\right\Vert +\frac{1}{4n}%
\sum_{i=1}^{n}\left( \int_{a}^{b}\frac{\left[ M_{i}\left( t\right)
-m_{i}\left( t\right) \right] ^{2}}{M_{i}\left( t\right) +m_{i}\left(
t\right) }dt\right) .  \label{3.18}
\end{equation}%
The equality holds in (\ref{3.18}) if and only if%
\begin{equation*}
\int_{a}^{b}\left\Vert f\left( t\right) \right\Vert dt\geq \frac{1}{4n}%
\sum_{i=1}^{n}\left( \int_{a}^{b}\frac{\left[ M_{i}\left( t\right)
-m_{i}\left( t\right) \right] ^{2}}{M_{i}\left( t\right) +m_{i}\left(
t\right) }dt\right) 
\end{equation*}%
and%
\begin{equation*}
\int_{a}^{b}f\left( t\right) dt=\left( \int_{a}^{b}\left\Vert f\left(
t\right) \right\Vert dt-\frac{1}{4n}\sum_{i=1}^{n}\left( \int_{a}^{b}\frac{%
\left[ M_{i}\left( t\right) -m_{i}\left( t\right) \right] ^{2}}{M_{i}\left(
t\right) +m_{i}\left( t\right) }dt\right) \right) \sum_{i=1}^{n}e_{i}.
\end{equation*}
\end{corollary}

\begin{proof}
As in the proof of Corollary \ref{c2.5}, (\ref{3.17}), implies that%
\begin{equation*}
\left\Vert f\left( t\right) \right\Vert -\func{Re}\left\langle f\left(
t\right) ,e_{i}\right\rangle \leq \frac{1}{4}\cdot \frac{\left[ M_{i}\left(
t\right) -m_{i}\left( t\right) \right] ^{2}}{M_{i}\left( t\right)
+m_{i}\left( t\right) }\ 
\end{equation*}%
for a.e. $t\in \left[ a,b\right] $ and $i\in \left\{ 1,\dots ,n\right\} .$

Applying Theorem \ref{t3.1} for%
\begin{equation*}
M_{i}\left( t\right) :=\frac{1}{4}\cdot \frac{\left[ M_{i}\left( t\right)
-m_{i}\left( t\right) \right] ^{2}}{M_{i}\left( t\right) +m_{i}\left(
t\right) },\text{ \ \ }t\in \left[ a,b\right] \text{, \ }i\in \left\{
1,\dots ,n\right\} ,
\end{equation*}%
we deduce the desired result.
\end{proof}

\section{Applications for Complex-Valued Functions}

Let $e=\alpha +i\beta $ $\left( \alpha ,\beta \in \mathbb{R}\right) $ be a
complex number with the property that $\left| e\right| =1,$ i.e., $\alpha
^{2}+\beta ^{2}=1.$ The following proposition concerning a reverse of the
continuous triangle inequality for complex-valued functions may be stated:

\begin{proposition}
\label{p.4.1} Let $f:\left[ a,b\right] \rightarrow \mathbb{C}$ be a Lebesgue
integrable function with the property that there exists a constant $\rho \in
\left( 0,1\right) $ such that 
\begin{equation}
\left| f\left( t\right) -e\right| \leq \rho \text{ for a.e. }t\in \left[ a,b%
\right] ,  \label{e.4.1}
\end{equation}
where $e$ has been defined above. Then we have the following reverse of the
continuous triangle inequality 
\begin{eqnarray}
&&\left( 0\leq \right) \int_{a}^{b}\left| f\left( t\right) \right| dt-\left|
\int_{a}^{b}f\left( t\right) dt\right|  \label{e.4.2} \\
&\leq &\frac{\rho ^{2}}{\sqrt{1-\rho ^{2}}\left( 1+\sqrt{1-\rho ^{2}}\right) 
}\left[ \alpha \int_{a}^{b}\func{Re}f\left( t\right) dt+\beta \int_{a}^{b}%
\func{Im}f\left( t\right) dt\right] .  \notag
\end{eqnarray}
\end{proposition}

The proof follows by Corollary \ref{c2.2}, and the details are omitted.

On the other hand, the following result is perhaps more useful for
applications:

\begin{proposition}
\label{p.4.2} Assume that $f$ and $e$ are as in Proposition \ref{p.4.1}. If
there exists the constants $M\geq m>0$ such that either 
\begin{equation}
\func{Re}\left[ \left( Me-f\left( t\right) \right) \left( \overline{f\left(
t\right) }-m\overline{e}\right) \right] \geq 0  \label{e.4.3}
\end{equation}
or, equivalently, 
\begin{equation}
\left| f\left( t\right) -\frac{M+m}{2}e\right| \leq \frac{1}{2}\left(
M-m\right)  \label{e.4.4}
\end{equation}
for a.e. $t\in \left[ a,b\right] ,$ holds, then 
\begin{eqnarray}
&&\left( 0\leq \right) \int_{a}^{b}\left| f\left( t\right) \right| dt-\left|
\int_{a}^{b}f\left( t\right) dt\right|  \label{e.4.5} \\
&\leq &\frac{\left( \sqrt{M}-\sqrt{m}\right) ^{2}}{2\sqrt{Mm}}\left[ \alpha
\int_{a}^{b}\func{Re}f\left( t\right) dt+\beta \int_{a}^{b}\func{Im}f\left(
t\right) dt\right] .  \notag
\end{eqnarray}
\end{proposition}

The proof may be done on utilising Corollary \ref{c2.3}, but we omit the
details

\begin{remark}
\label{r.4.1} From a practical view point, since 
\begin{eqnarray*}
&&\func{Re}\left[ \left( Me-f\left( t\right) \right) \left( \overline{%
f\left( t\right) }-m\overline{e}\right) \right] \\
&=&\left[ M\alpha -\func{Re}f\left( t\right) \right] \left[ \func{Re}f\left(
t\right) -m\alpha \right] +\left[ M\beta -\func{Im}f\left( t\right) \right] %
\left[ \func{Im}f\left( t\right) -m\beta \right] ,
\end{eqnarray*}%
hence a sufficient condition for \ref{e.4.3} to hold is 
\begin{equation}
m\alpha \leq \func{Re}f\left( t\right) \leq M\alpha \text{ and }m\beta \leq 
\func{Im}f\left( t\right) \leq M\beta  \label{e.4.6}
\end{equation}%
for a.e. $t\in \left[ a,b\right] ,$ where $\alpha ,\beta $ are assumed to be
positive and satisfying the condition $\alpha ^{2}+\beta ^{2}=1$. We observe
that the above condition (\ref{e.4.6}) is very easy to verify in practice,
therefore it may useful in various applications where reverses of the
continuous triangle inequality are required.
\end{remark}

Finally, on making use of Corollary \ref{c2.5}, one may state the following
result as well:

\begin{proposition}
\label{p.4.3} Let $f$ be as in Proposition \ref{p.4.1} and the measurable
functions $K,k:\left[ a,b\right] \rightarrow \lbrack 0,\infty )$ with the
property that 
\begin{equation*}
\frac{\left( K-k\right) ^{2}}{K+k}\in L\left[ a,b\right]
\end{equation*}%
and 
\begin{equation*}
\alpha k\left( t\right) \leq \func{Re}f\left( t\right) \leq \alpha K\left(
t\right) \text{ and }\beta k\left( t\right) \leq \func{Im}f\left( t\right)
\leq \beta K\left( t\right)
\end{equation*}%
for a.e. $t\in \left[ a,b\right] ,$ where $\alpha ,\beta $ are assumed to be
positive and satisfying the condition $\alpha ^{2}+\beta ^{2}=1$. Then the
following reverse of the continuous triangle inequality is valid: 
\begin{eqnarray*}
&&\left( 0\leq \right) \int_{a}^{b}\left\vert f\left( t\right) \right\vert
dt-\left\vert \int_{a}^{b}f\left( t\right) dt\right\vert \\
&\leq &\frac{1}{4}\int_{a}^{b}\frac{\left[ K\left( t\right) -k\left(
t\right) \right] ^{2}}{K\left( t\right) +k\left( t\right) }dt.
\end{eqnarray*}%
The constant $\frac{1}{4}$ is best possible in the sense that it cannot be
replaced by a smaller quantity.
\end{proposition}

\begin{remark}
\label{r.4.2} One may realise that similar results can be stated if the
Corollaries \ref{c3.2}-\ref{c3.5}\ stated above are used. For the sake of
brevity, we do not state them here.
\end{remark}


\begin{thebibliography}{9}
\bibitem{DM} J.B. DIAZ and F.T. METCALF, A complementary triangle inequality
in Hilbert and Banach spaces, \textit{Proceedings Amer. Math. Soc., }\textbf{%
17}(1) (1966), 88-97.

\bibitem{SSD1} S.S. DRAGOMIR, Reverses of the continuous triangle inequality
for Bochner integral of vector valued function in Hilbert spaces. \textit{%
RGMIA Res. Rep. Coll., to appear.}

\bibitem{K} J. KARAMATA, \textit{Teorija i Praksa Stieltjesova Integrala}
(Serbo-Coratian) (Stieltjes Integral, Theory and Practice), SANU, Posebna
izdanja, 154, Beograd, 1949.

\bibitem{M} M. MARDEN, The Geometry of the Zeros of a Polynomial in a
Complex Variable, \textit{Amer. Math. Soc. Math. Surveys}, \textbf{3}, New
York, 1949.

\bibitem{MPF} D.S. MITRINOVI\'{C}, J.E. PE\v{C}ARI\'{C} and\ A.M. FINK, 
\textit{Classical and New Inequalities in Analysis, }Kluwer Academic
Publishers, Dordrecht/Boston/London, 1993.

\bibitem{P} M. PETROVICH, Module d'une somme, \textit{L' Ensignement Math%
\'{e}matique,} \textbf{19} (1917), 53-56.

\bibitem{W} H.S. WILF, Some applications of the inequality of arithmetic and
geometric means to polynomial equations, \textit{Proceedings Amer. Math.
Soc., }\textbf{14} (1963), 263-265.
\end{thebibliography}
\end{document}